\documentclass[12pt,a4paper]{article}

\usepackage{fancyhdr, amsmath, amssymb, amsfonts, dsfont, subfigure,amsthm,tikz,color,epsfig,graphicx,soul,mathtools}

\usepackage{IEEEtrantools}

\usetikzlibrary{calc,backgrounds}

\newlength{\guillotine}
	\newtheorem{thm}{Theorem}
	\newtheorem{cor}[thm]{Corollary}

	\newtheorem{prop}[thm]{Proposition}

	\newtheorem{defn}[thm]{Definition}

	\theoremstyle{remark}

	\newcommand\RR{\mathbb R}

	\newcommand\calj{\mathcal J}
	\newcommand\calg{\mathcal G}
	
	\newcommand\cali{\mathcal I}

    \newcommand\alf\alpha
	\newcommand\de\delta
	\newcommand\eps\varepsilon
	\newcommand\la\lambda
	\newcommand\De\Delta
	\newcommand\La\Lambda
	\newcommand\Om\Omega
	\newcommand\kap\kappa
    \newcommand\sig\sigma

	\makeatletter
	\newcommand{\ostar}{\mathbin{\mathpalette\make@circled\star}}
	\newcommand{\make@circled}[2]{%
	  \ooalign{$\m@th#1\smallbigcirc{#1}$\cr\hidewidth$\m@th#1#2$\hidewidth\cr}%
	}
	\newcommand{\smallbigcirc}[1]{%
	  \vcenter{\hbox{\scalebox{0.77778}{$\m@th#1\bigcirc$}}}%
	}
	\makeatother

	\newcommand\emp\emptyset

    \newcommand\un\underline
    
    \newcommand\ol\overline

	\DeclareMathOperator\per{per}
	\DeclareMathOperator\In{in}
	\DeclareMathOperator\vol{vol}
	\DeclareMathOperator\Int{int}
	\DeclareMathOperator\nab\nabla

\title{Simple bounds for the inradius and inner neighbourhood of a convex body}
\author{%
	Benedict Sewell%
	\thanks{Supported by the Alfr\'ed R\'enyi young researcher fund.}
}
\date{\today}

\begin{document}

\maketitle

\begin{abstract}
	In this short note, we show that the inradius of a convex body in $\RR^n$ is comparable to its volume divided by its surface area. We also give a simple formula, in terms of the volume and inradius, that is comparable to the volume of the intersection of a convex body with the $\eps$-neighbourhood of its boundary, and provide an application of this to self-projective sets with convex holes, extending results of \cite{sewell}.
\end{abstract}

We start with some simple definitions.
\begin{defn}
	Fix a convex body $\Om\subset \RR^n$  (a convex, closed set with non-empty interior); let
		\begin{itemize}

			\item $\vol(\Om)$ denote its $n$-dimensional Lebesgue measure;

			\item $\per(\Om)$ denote the $(n-1)$-dimensional Lebesgue measure of its boundary, $\partial \Om$; and

			\item $\In(\Om)$ denote its \emph{inradius}:

			\subitem--- the radius of the largest inscribed ball within $\Om$, or equivalently

			\subitem--- the furthest distance from the boundary of any point in $\Om$.

		\end{itemize}
	We call a point in $\Om$ maximally distant from $\partial\Om$ an \emph{incentre} of $\Om$ and \emph{the} incentre if it is unique.
\end{defn}

In this short note, we mainly concern ourselves with the volume of the $\eps$-\emph{inner neighbourhood} $L_\eps(\Om)$, defined by
	 $$
		 		L_\eps(\Om)
		 	=
		 		\big\{
		 			x \in \Om
		 		:
		 			\|x - y\|
		 		\leq
		 			\eps
		 		\text{ for some }
		 			y \in \partial\Om
		 		\big\}
		 	;
	 $$
that is, the set of points in $\Om$ at most $\eps$ away from the boundary $\partial \Om$.

Such a formula will feature the volume of $\Om$ and its inradius, and so it appears expedient that we should relate the latter to other simple quantities. This is the object of the first section, which also allows us to introduce some of the tools used.

The second section shows that $g(\eps) \geq \vol(L_\eps(\Om)) \geq g(\eps)/n$, where
	$$
			g(\eps) = g(\eps,\Om)
		:=
			\vol(\Om)
			\left(
				1
			-
				\max
				\left(
				0,
					1
				-
					\frac
						{\eps}
						{\In(\Om)}
				\right)^n
			\right).
	$$

In the third and final section, we state extensions of this inequality to the main results of \cite{sewell} which follow from the main formula.

\section{The inradius}

In this section we state and prove the formula below, generalising Heron's formula for the inradius of a triangle.

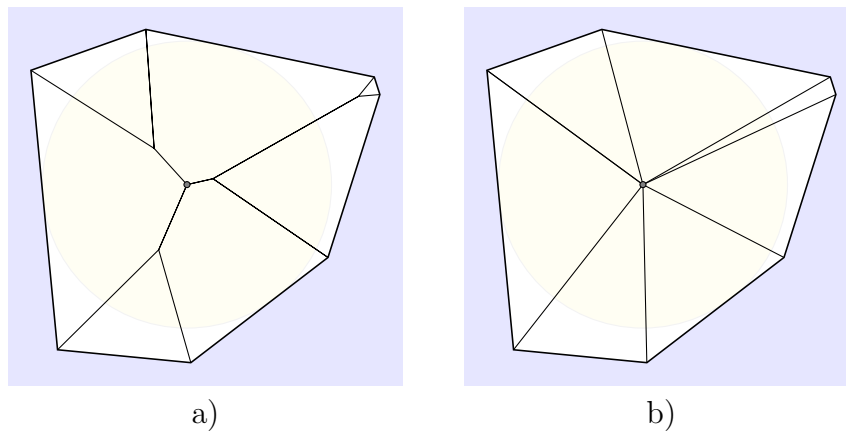
\begin{figure}[ht]

	\centering
	\begin{tikzpicture}[scale = 6]

		\path[fill = blue!10] (-0.00616, 0.08352) rectangle (0.85902,0.91987);

		\path[fill = blue!10] (-0.00616, 0.08352) ++ (1,0) rectangle (1.85902,0.91987);

		\draw[semithick,fill = white] (0.04384, 0.77957) --
				(0.10215, 0.16202) --
				(0.39443, 0.13352) --
				(0.69521, 0.36598) --
				(0.80902, 0.72602) --
				(0.79672, 0.76476) --
				(0.29577, 0.86987) --
				(0.04384, 0.77957) --
				cycle;

		\node (B) at (0.386, 0.527) [circle,fill=black,inner sep = 0.5] {};

		\draw[fill = yellow, opacity = 0.05] (B) circle[radius = 0.317];

		\draw[semithick, fill = white] (1.04384, 0.77957) --
				(1.10215, 0.16202) --
				(1.39443, 0.13352) --
				(1.69521, 0.36598) --
				(1.80902, 0.72602) --
				(1.79672, 0.76476) --
				(1.29577, 0.86987) --
				(1.04384, 0.77957) --
				cycle;

		\node (A) at (1.386, 0.527) [circle,fill=black,inner sep = 0.5] {};

		\draw[fill = yellow, opacity = 0.05] (A) circle[radius = 0.317];

		\foreach \x in {{(0.04384, 0.77957)},{(0.10215, 0.16202)},{(0.39443, 0.13352)},{(0.69521, 0.36598)},{(0.80902, 0.72602)},{(0.79672, 0.76476)},{(0.29577, 0.86987)},{(0.04384, 0.77957)}}
		{
			\draw \x ++ (1,0) -- (A);
		}

	 	\draw (0.10215, 0.16202) -- (0.324, 0.383) -- (B) -- (0.324, 0.383) -- (0.39443, 0.13352);
	 	% lower left triangle

	 	\draw (0.80902, 0.72602) -- (0.762, 0.722) -- (0.443, 0.54) -- (B) -- (0.443, 0.54) -- (0.69521, 0.36598) -- (0.443, 0.54) -- (0.762, 0.722) -- (0.79672, 0.76476);
	 	% right part conglomeration

	 	\draw (0.04384, 0.77957) -- (0.314, 0.607) -- (0.29577, 0.86987) -- (0.314, 0.607) -- (B);

		\node at (A) [circle, draw = black, fill = gray, inner sep = 0.8] {};

		\node at (B) [circle, draw = black, fill = gray, inner sep = 0.8] {};

		\path (-0.00616, 0.08352) -- (0.85902,0.08352) node[midway, anchor = north, inner sep = 5pt] {a)};

		\path (-0.00616, 0.08352) ++ (1,0) -- (1.85902,0.08352) node[midway, anchor = north, inner sep = 5pt] {b)};

	\end{tikzpicture}
	\caption{Proof of Proposition \ref{prop-heron}: the decomposition of $\Om$ into a) $\{\Om_S\}$ and b) $\{\De_S\}$. The largest inscribed ball is depicted in faint yellow, and the incentre with a dot.}
	\label{fig:1}
\end{figure}

\begin{prop}
	\label{prop-heron}
	For any convex body $\Om\subset\RR^n$,
		$$
%				\frac1n
				\frac
					{\vol(\Om)}
					{\per(\Om)}
			\leq
				\In(\Om)
			\le
				n
				\frac
					{\vol(\Om)}
					{\per(\Om)}.
		$$
	Moreover, the constants cannot be improved.
\end{prop}

	Since a convex body can be approximated arbitrarily well in the Hausdorff metric by convex polytopes, it suffices to prove this result when $\Om$ is itself a polytope. The proof then relies on two useful decompositions for $\Om$, which are illustrated in Figure \ref{fig:1}. These will be used in later proofs as well, and may naturally correspond.

\begin{proof}[Proof of lower bound]
		Consider one side $S\subset \partial\Om$ of $\Om$, and the locus $\Om_S\subset \Om$ of points whose closest side of $\Om$ is $S$, as in Figure 1a). Treating $S \subset \RR^{n-1}\times \{0\}$ and $\Om \subset \RR^{n-1}\times [0,\infty)$, we see that
			\begin{itemize}

				\item every point in $\Om_S$ has distance at most $\In(\Om)$ distance from $S$, and

				\item $\Om_S$ is bound by the hyperplanes bisecting $S$ and one of its adjacent sides:

				by convexity, the interior angle between two sides lies in $(0,\pi)$, and thus the bisector makes an internal angle less than $\pi/2$ with $S$: i.e., is sloping inwards towards $S$.\footnote{We will give a more precise description later.}
			\end{itemize}
		%
%		(We will give a more precise description later.)
		Hence $\Om_S \subset S \times [0,\In(\Om)]$, in such a way that
			$$
					\vol(\Om_S)
				<
					\vol_{n-1}(S)\cdot\In(\Om)
			$$
		holds. Summing over $S$ gives the required inequality.
\end{proof}
\begin{proof}[Proof of upper bound]
		First fix an incentre $x\in \Om$.
		Again consider the side $S$, and let $\De_S$ denote the cone with base $S$ and vertex $x$.
		Again taking $S \subset \RR^{n-1}\times \{0\}$ and $\Om \subset \RR^{n-1}\times [0,\infty)$ for simplicity;
		since the ball of radius $\In(\Om)$ based at $x$ is contained in $\Om \subset \RR^{n-1}\times [0,\infty)$,
		we must have $x \in \RR^{n-1}\times [\In(\Om),\infty)$.
		Therefore, since $\De_S$ is a cone with base $S$ and height at least $\In(\Om)$,
			\begin{align}
					\frac1n
					\vol_{n-1}(S)\cdot\In(\Om)
				\leq
					\vol(\De_S)
				,
			\label{eq:delta-S-decomp}
			\end{align}
		with equality when $x\in \RR^{n-1} \times \{\In(\Om)\}$, i.e., when the maximal ball based at $x$ meets $\RR^{n-1}\times\{0\}$, and hence $S$ in order to avoid a contradiction.
		The proof is again finished by  summing over $S$ in \eqref{eq:delta-S-decomp}.
\end{proof}

Finally, to conclude the proof, we present a simple family of examples.
\begin{proof}[Proof that constants are optimal]
	Consider the ``square pancake" $\Om  = [0,1] \times [0,K]^{n-1}$ for $n\geq 2$ and $K\geq 1$:
		\begin{itemize}
			\item $\In(\Om) = \frac12$.

			\item $\vol(\Om) = K^{n-1}$.

			\item $\per(\Om) = 2 \big( K^{n-1} + (n-2) K^{n-2} + 1\big)$.
		\end{itemize}
	I.e.,
		$$
				\frac
					{\vol(\Om)}
					{\per(\Om)}
			=
				\frac12
				\frac
					{K^{n-1}}
					{K^{n-1} + (n-2) K^{n-2} + 1}
			=
				\frac
					{\In(\Om)}
					{1 + (n-2) K^{-1} + K^{1-n}}
			.
		$$
	The right hand side equals $\In(\Om)/n$ when $K=1$, and approaches $\In(\Om)$ as $K \to \infty$.
\end{proof}

This completes the proof of the proposition. We end this section by stating what polytopes attain the upper bound, following immediately from the above proof.

The following is a special case of a result of Aspergen \cite{aspegren}, which was originally shown for $n=2$ by Apostol and Mnatsakanian \cite{apostol}.

\begin{cor}[\cite{aspegren}] The convex polytope $\Om$ satisfies
		$$
				\In(\Om)
			=
				n
			\cdot
				\frac
					{\vol(\Om)}
					{\per(\Om)}
		$$
	if {and only if} it is \emph{circumscribed}: there exists a ball contained within $\Om$ that meets all sides.\footnote{It follows in particular that it is the unique, maximal ball in $\Om$.}
\end{cor}
	To be precise, Aspegren showed in \cite{aspegren} that this formula applies to the more general class of all convex sets for which every point on the boundary has a tangent plane tangent to the maximal ball.

\section{The inner neighbourhood}

\begin{figure}[htb]
	\centering
	\includegraphics{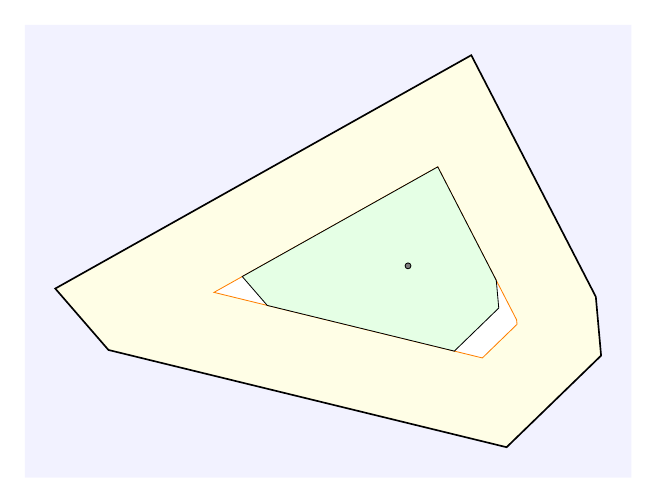}
	\caption{Illustration of Proposition \ref{prop:scale-copy}. The scale copy of $\Om$ is in green, and $L_\eps(\Om)$ in yellow, for a certain polygon $\Om$ and value of $\eps \in (0, \In(\Om))$. The incentre is marked with a dot.}
	\label{fig:scale-copy}
\end{figure}

We turn our attention to the volume of the $\eps$-inner neighbourhood, as defined previously and also depicted in Figure \ref{fig:scale-copy} in a certain case.

Our first result uses that the closure of the complement $\Om\setminus L_{\eps}(\Om)$ contains a large scale copy of $\Om$, as illustrated in Figure \ref{fig:scale-copy}.

\begin{prop}
	For any $\eps \in[0,\In(\Om)]$,
		\begin{equation}
				L_\eps(\Om)
			\subset
%				\begin{cases}
%					\Om,
%				&
%					\text{if }\eps \geq \In(\Om),
%				\\
					\Om
				\setminus
					\Int\left[
						\left(
							1-
							 \frac
							 	{\eps}
							 	{\In(\Om)}
						\right)
					\cdot
						\Om
					\right]
					\!\!,
%				&
%					\text{if }\eps < \In(\Om);
%				\end{cases}
		\label{eq:L_epsilon-made-explicit}
		\end{equation}
	where $x \mapsto \la\cdot x$ denotes enlargement, scale factor $\la$, about a fixed incentre of $\Om$.
	\label{prop:scale-copy}
\end{prop}

The proof of this result is a straightforward application of convexity:

\begin{proof}
	We first show $f : x \mapsto d(x,\partial\Om)$ is concave on $\Om$, as follows. Suppose that $B_0$, $B_1$ are the largest balls in $\Om$ based respectively at $x_0$, $x_1 \in \Om$ (i.e., with radii $f(x_0)$ and $f(x_1)$). Then, for any $\la\in(0,1)$,
		$$
			B_\la := \la B_1 + (1-\la) B_0
		$$
	is a ball based at $x_\la := \la x_1 + (1-\la) x_0$, with radius $\la f(x_1) + (1-\la)f(x_0)$, contained wholly in $\Om$ by convexity. Hence,
		$$
				f(x_\la)
			\geq
				\la f(x_1) + (1-\la)f(x_0),
		$$
	as claimed.
	Therefore, given $y \in \Int(\Om)$ and $\eps < f(y)$, contracting $x\in \Om$ towards $y$ by a factor $\left(1 - \frac\eps{f(y)}\right)$ takes it outside of the $\eps$-inner neighbourhood:
		$$
				f %d_{\partial\Om}
				\left(
					\left(
						\frac
							\eps
							{f(y)}%{d_{\partial\Om}(y)}
					\right)
						y
					+
					\left(
						1
					-
						\frac
							\eps
							{f(y)}
							%{d_{\partial\Om}(y)}
					\right)
					x
				\right)
			\geq
					\left(
						\frac
							\eps
							{f(y)}
							%{d_{\partial\Om}(y)}
					\right)
					f(y)
					%d_{\partial\Om}(y)
				+
					\left(
						1
					-
						\frac
							\eps
							{f(y)}
							%{d_{\partial\Om}(y)}
					\right)
						f(x)
						%d_{\partial\Om}(x)
			\geq \eps.
		$$
	Thus taking $y$ to be an incentre, i.e., so that $f(y) = \In(\Om)$, gives \eqref{eq:L_epsilon-made-explicit convex verzio}.
\end{proof}

The important immediate consequence is the following inequality, which introduces the formula $g$ which implicitly depends on the volume and inradius of $\Om$.
\begin{prop}
	For any $\eps \in [0,\In(\Om)]$,
	\begin{align}
			\vol(L_\eps(\Om))
		&\leq
			g(\eps)
		:=
			\vol(\Om)
			\left(
				1
			-
				\left(
						1
					-
						\frac
							\eps
							{\In(\Om)}
					\right)^{\!\!n}
				\right)
				.
	\label{eq:L_eps-upper bound}
	\end{align}
\end{prop}

	Before we continue, we note that circumscribed polytopes (and probably their above generalisation involving tangent planes) again play a special role in giving equality in \eqref{eq:L_epsilon-made-explicit} and \eqref{eq:L_eps-upper bound}.

\begin{prop}
	When $\Om$ is a circumscribed polytope, $\vol(L_\eps(\Om)) = g(\eps)$ for all $\eps\in[0,\In(\Om)]$.
\end{prop}

\begin{proof}
	Given a side $S\subset\Om$, we first claim that $\Om_S = \De_S$ is a cone with base $S$ and vertex at the unique incentre. This follows simply because all the hyperplanes bisecting $S$ and its adjacent sides necessarily meet the incentre.

	In particular, the enlargement spoken of above maps $\Om_S$ into itself, and is more precisely an enlargement towards its vertex, the incentre.

	In particular, if $x \in L_\eps(\Om)$, treating $S$ as a flat something.
\end{proof}

	We note again that circumscribed polytopes play a special role in giving equality in the above, as can be shown by noticing that the enlargement takes the locus $\Om_S=\De_S$, a cone with base $S$, into itself.

For our application to projective fractal geometry, it is sufficient to show that the two sides are equal up to a constant. Indeed, as hinted above, one obtains equality for

We obtain such a bound simply via a linear one, following from the concavity of $\eps \to \vol(L_\eps(\Om)).$

\begin{prop}
	When $\Om$ is a polytope, the function $\eps \mapsto \vol(L_\eps(\Om))$ is continuously differentiable with a decreasing derivative. Subsequently, it is concave and
		\begin{align}
				\vol(L_\eps(\Om))
			&\geq
				\eps
%			\cdot
		\,
				\frac{\vol(\Om)}{\In(\Om)}
			\geq
				\frac
					{g(\eps)}
					n
%				\frac
%					{\vol(\Om)}
%					n
%				\left(
%					1
%				-
%					\left(
%						1
%					-
%						\frac
%							\eps
%							{\In(\Om)}
%					\right)^{\!\! n}
%				\right)
		\end{align}
	for all $\eps \in [0,\In(\Om)]$.
	\label{prop-curves}
\end{prop}
\begin{figure}
	\centering
	\begin{tikzpicture}[yscale=3,xscale=5]

%				\fill[yellow!50] (1.5,1) -- (1,1) parabola (0,0) -- (1.5,0);
				\fill[yellow!50] (0.,0.)--(0.01,0.029701)--(0.02,0.058808)--(0.03,0.087327)--(0.04,0.115264)--(0.05,0.142625)--(0.06,0.169416)--(0.07,0.195643)--(0.08,0.221312)--(0.09,0.246429)--(0.1,0.271)--(0.11,0.295031)--(0.12,0.318528)--(0.13,0.341497)--(0.14,0.363944)--(0.15,0.385875)--(0.16,0.407296)--(0.17,0.428213)--(0.18,0.448632)--(0.19,0.468559)--(0.2,0.488)--(0.21,0.506961)--(0.22,0.525448)--(0.23,0.543467)--(0.24,0.561024)--(0.25,0.578125)--(0.26,0.594776)--(0.27,0.610983)--(0.28,0.626752)--(0.29,0.642089)--(0.3,0.657)--(0.31,0.671491)--(0.32,0.685568)--(0.33,0.699237)--(0.34,0.712504)--(0.35,0.725375)--(0.36,0.737856)--(0.37,0.749953)--(0.38,0.761672)--(0.39,0.773019)--(0.4,0.784)--(0.41,0.794621)--(0.42,0.804888)--(0.43,0.814807)--(0.44,0.824384)--(0.45,0.833625)--(0.46,0.842536)--(0.47,0.851123)--(0.48,0.859392)--(0.49,0.867349)--(0.5,0.875)--(0.51,0.882351)--(0.52,0.889408)--(0.53,0.896177)--(0.54,0.902664)--(0.55,0.908875)--(0.56,0.914816)--(0.57,0.920493)--(0.58,0.925912)--(0.59,0.931079)--(0.6,0.936)--(0.61,0.940681)--(0.62,0.945128)--(0.63,0.949347)--(0.64,0.953344)--(0.65,0.957125)--(0.66,0.960696)--(0.67,0.964063)--(0.68,0.967232)--(0.69,0.970209)--(0.7,0.973)--(0.71,0.975611)--(0.72,0.978048)--(0.73,0.980317)--(0.74,0.982424)--(0.75,0.984375)--(0.76,0.986176)--(0.77,0.987833)--(0.78,0.989352)--(0.79,0.990739)--(0.8,0.992)--(0.81,0.993141)--(0.82,0.994168)--(0.83,0.995087)--(0.84,0.995904)--(0.85,0.996625)--(0.86,0.997256)--(0.87,0.997803)--(0.88,0.998272)--(0.89,0.998669)--(0.9,0.999)--(0.91,0.999271)--(0.92,0.999488)--(0.93,0.999657)--(0.94,0.999784)--(0.95,0.999875)--(0.96,0.999936)--(0.97,0.999973)--(0.98,0.999992)--(0.99,0.999999)--(1.,1.) -- (0,0);

				\draw[dashed, gray] (0,0) rectangle (1,1);

				\draw[thick]
					(1.5,1)
				-- 	(1,1)
				-- (0,0);

				\draw[thick] (0.,0.)--(0.01,0.00990033)--(0.02,0.0196027)--(0.03,0.029109)--(0.04,0.0384213)--(0.05,0.0475417)--(0.06,0.056472)--(0.07,0.0652143)--(0.08,0.0737707)--(0.09,0.082143)--(0.1,0.0903333)--(0.11,0.0983437)--(0.12,0.106176)--(0.13,0.113832)--(0.14,0.121315)--(0.15,0.128625)--(0.16,0.135765)--(0.17,0.142738)--(0.18,0.149544)--(0.19,0.156186)--(0.2,0.162667)--(0.21,0.168987)--(0.22,0.175149)--(0.23,0.181156)--(0.24,0.187008)--(0.25,0.192708)--(0.26,0.198259)--(0.27,0.203661)--(0.28,0.208917)--(0.29,0.21403)--(0.3,0.219)--(0.31,0.22383)--(0.32,0.228523)--(0.33,0.233079)--(0.34,0.237501)--(0.35,0.241792)--(0.36,0.245952)--(0.37,0.249984)--(0.38,0.253891)--(0.39,0.257673)--(0.4,0.261333)--(0.41,0.264874)--(0.42,0.268296)--(0.43,0.271602)--(0.44,0.274795)--(0.45,0.277875)--(0.46,0.280845)--(0.47,0.283708)--(0.48,0.286464)--(0.49,0.289116)--(0.5,0.291667)--(0.51,0.294117)--(0.52,0.296469)--(0.53,0.298726)--(0.54,0.300888)--(0.55,0.302958)--(0.56,0.304939)--(0.57,0.306831)--(0.58,0.308637)--(0.59,0.31036)--(0.6,0.312)--(0.61,0.31356)--(0.62,0.315043)--(0.63,0.316449)--(0.64,0.317781)--(0.65,0.319042)--(0.66,0.320232)--(0.67,0.321354)--(0.68,0.322411)--(0.69,0.323403)--(0.7,0.324333)--(0.71,0.325204)--(0.72,0.326016)--(0.73,0.326772)--(0.74,0.327475)--(0.75,0.328125)--(0.76,0.328725)--(0.77,0.329278)--(0.78,0.329784)--(0.79,0.330246)--(0.8,0.330667)--(0.81,0.331047)--(0.82,0.331389)--(0.83,0.331696)--(0.84,0.331968)--(0.85,0.332208)--(0.86,0.332419)--(0.87,0.332601)--(0.88,0.332757)--(0.89,0.33289)--(0.9,0.333)--(0.91,0.33309)--(0.92,0.333163)--(0.93,0.333219)--(0.94,0.333261)--(0.95,0.333292)--(0.96,0.333312)--(0.97,0.333324)--(0.98,0.333331)--(0.99,0.333333)--(1.,0.333333)--(1.5,1/3) node[anchor = west] {$g/n$}; % bottom curve

				\draw[<->] (0,1.7) -- (0,0) -- (1.7,0); %node[anchor = north east]{$\eps$};

				\draw[thick] (0.,0.)--(0.01,0.029701)--(0.02,0.058808)--(0.03,0.087327)--(0.04,0.115264)--(0.05,0.142625)--(0.06,0.169416)--(0.07,0.195643)--(0.08,0.221312)--(0.09,0.246429)--(0.1,0.271)--(0.11,0.295031)--(0.12,0.318528)--(0.13,0.341497)--(0.14,0.363944)--(0.15,0.385875)--(0.16,0.407296)--(0.17,0.428213)--(0.18,0.448632)--(0.19,0.468559)--(0.2,0.488)--(0.21,0.506961)--(0.22,0.525448)--(0.23,0.543467)--(0.24,0.561024)--(0.25,0.578125)--(0.26,0.594776)--(0.27,0.610983)--(0.28,0.626752)--(0.29,0.642089)--(0.3,0.657)--(0.31,0.671491)--(0.32,0.685568)--(0.33,0.699237)--(0.34,0.712504)--(0.35,0.725375)--(0.36,0.737856)--(0.37,0.749953)--(0.38,0.761672)--(0.39,0.773019)--(0.4,0.784)--(0.41,0.794621)--(0.42,0.804888)--(0.43,0.814807)--(0.44,0.824384)--(0.45,0.833625)--(0.46,0.842536)--(0.47,0.851123)--(0.48,0.859392)--(0.49,0.867349)--(0.5,0.875)--(0.51,0.882351)--(0.52,0.889408)--(0.53,0.896177)--(0.54,0.902664)--(0.55,0.908875)--(0.56,0.914816)--(0.57,0.920493)--(0.58,0.925912)--(0.59,0.931079)--(0.6,0.936)--(0.61,0.940681)--(0.62,0.945128)--(0.63,0.949347)--(0.64,0.953344)--(0.65,0.957125)--(0.66,0.960696)--(0.67,0.964063)--(0.68,0.967232)--(0.69,0.970209)--(0.7,0.973)--(0.71,0.975611)--(0.72,0.978048)--(0.73,0.980317)--(0.74,0.982424)--(0.75,0.984375)--(0.76,0.986176)--(0.77,0.987833)--(0.78,0.989352)--(0.79,0.990739)--(0.8,0.992)--(0.81,0.993141)--(0.82,0.994168)--(0.83,0.995087)--(0.84,0.995904)--(0.85,0.996625)--(0.86,0.997256)--(0.87,0.997803)--(0.88,0.998272)--(0.89,0.998669)--(0.9,0.999)--(0.91,0.999271)--(0.92,0.999488)--(0.93,0.999657)--(0.94,0.999784)--(0.95,0.999875)--(0.96,0.999936)--(0.97,0.999973)--(0.98,0.999992)--(0.99,0.999999)--(1.,1.);

				\node at (1,-0.05) [anchor = north] {$\In(\Om)$};

				\node at (1.5,1) [anchor = west] {$g$};

				\node at (0,1) [anchor = east] {$\mathllap{\vol(\Om)}$};

				\begin{pgfonlayer}{background}
					\fill[blue!1]
						(current bounding box.south west) ++ (-0.5,-0.05) rectangle
						(current bounding box.north east) ++ (0.5,0.05)  rectangle (current bounding box.south west);
				\end{pgfonlayer}

	\end{tikzpicture}
	\caption{Sketch of Proposition \ref{prop-curves}: $\vol(L_\eps(\Om))$ is constrained to lie in the yellow region.}
\end{figure}
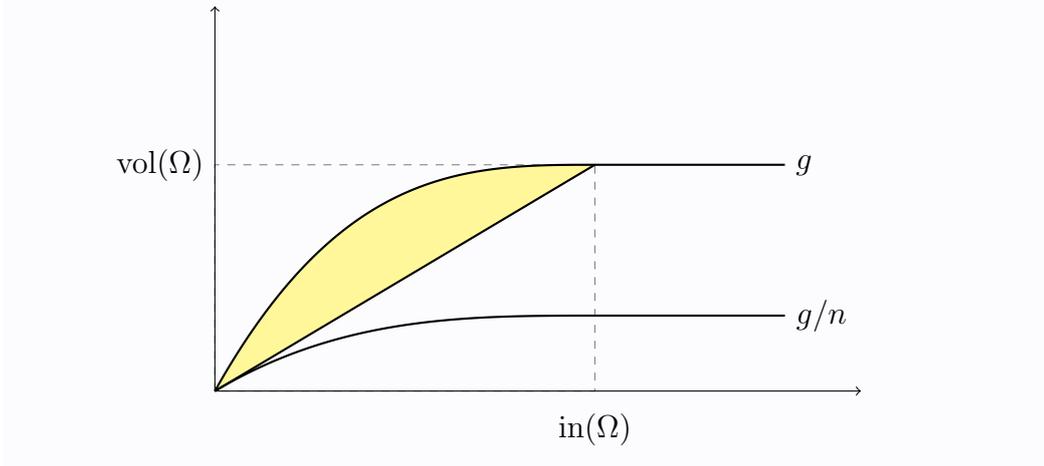

\begin{proof}
		We prove the statements in reverse order. First, note that the rightermost inequality is simply a rearrangement of the Bernoulli inequality, and is not related to functional properties of $\eps\mapsto \vol(L_\eps(\Om))$. The other inequality follows from concavity, since we know that the graph passes through $(0,0)$ and $(\In(\Om),\vol(\Om))$, and is hence bounded below by the line joining these two.

		To show the claimed functional properties, we again consider the case that $\Om$ is a polytope, and recall the locus $\Om_S\subset \Om$ of those points whose closest side is $S$. Furthermore, let $\Pi_S$ denote its supporting hyperplane (and similarly for other sides) We denote $H_{S,S'}$ as the half space comprising those points either closer to $S$ than $S'$, or equidistant.

		We give an alternative description of $\Om_S$ as the set of those points in $\Om$ whose closest hyperplane (out of those supporting sides) is $\Pi_S$:
			$$
					\Om_S
				=
						\Om
					\cap
						\bigcap_{S'\neq S}
							H_{S'\geq S}.
			$$
		This is easy to show, since for any given point in $\Om$, the closest point on a supporting hyperplane cannot lie outside $\Om$, otherwise there would be an intermediate point on the boundary $\partial \Om$, hence on another supporting hyperplane.

		In particular, this formula shows that $\Om_{S'}$ is a convex polytope, and we will also use it to show that the $(n-1)$-dimensional volume of the set of those points in $\Om_S$ distance $\eps$ from $S$ (equivalently $\partial\Om$) is both continuous and decreasing in $\eps$. This is sufficient to prove the proposition, since this volume equals
			$$
				\frac
					{\mathrm d}
					{\mathrm d\eps}
				\vol(L_\eps(\Om)\cap\Om_S).
			$$
		and the required functional properties follow by summation.

		Continuity follows naturally, and the decreasing property by induction:

		Taking the intersection of any plane $\Pi$ orthogonal to $S$, we see that the intersection of $\Om_S\cap\Pi$ has the same description in one dimension lower: namely, since the bisector of $S$ and another side cannot be vertical, $H_{S,S'}\cap \Pi$ is a half-space within $\Pi$, and hence
			$$
					\Om_S \cap \Pi
				=
					\Om \cap \Pi
					\cap
					\bigcap_{S'\neq S}
						H_{S,S'} \cap \Pi
				.
			$$
		This allows us to inductively reduce the dimension of the problem until the base case, $n=2$; which is rather natural, since the description of the boundary of $\Om_S$ is essentially comprised of two piecewise-linear decreasing curves.
%
%		\begin{tikzpicture}
%			\draw (0,0) --
%		\end{tikzpicture}
%%
%
%
%
%
%
%%and we use the $\Om_S$ partition for a given side $S$ of $\Om$. Again treating $S\subset \RR^{n-1}\times\{0\}$ and $\Om_S \subset \RR\times[0,\In(\Om)]$, its contribution to $\vol(L_\eps)$ is simply $\vol(\Om_S \cap [0,\eps])$. In particular, it is differentiable with derivative $\vol_{n-1}(\Om_S \cap \{\eps\})$, which we show is continuous and decreasing:
%
%		$\Om_S$ is a convex polytope obtained by intersecting $S\times [0,\In(\Om)]$ with halfspaces whose hyperplanes do not intersect the interior of $S$ (i.e., given by the bisectors of the planes corresponding to $S$ and every other side). By induction, the volume of the $\eps$-slice is decreasing in $\eps$, and is naturally continuous:
%
%		Indeed, taking the heights $0 = \eps_0 < \eps_1 < \cdots < \eps_m$ for which $\Om_S$ has a vertex, the subset $\Om_S \cap[\eps_{k-1},\eps_{k}]$ has all the vertices on its top or bottom face, it is easy to show that the slice $H_y = \Om_S\cap\{y\}$ at height $y = \la \eps_{k-1} + (1-\la) \eps_k$ is precisely the corresponding Minkowski interpolation of its top and bottom:
%		%
%			 $$
%					H_y
%				=
%					\la
%					H_{\eps_{k-1}}
%				+
%					(1-\la)
%					H_{\eps_k}.
%			$$
%		%
%		Hence, by Minkowski's theorem, its $(n-1)$-dimensional volume is given by a degree-$(n-1)$ polynomial in $\la$, interpolating between $\vol(H_{\eps_{k-1}})$ and $\vol(H_{\eps_k})$.
\end{proof}

\section{Application to self-projective attractors}

Using the result of the previous section, the method of proof of the main theorem of \cite{sewell} naturally extends to obtain the following:

\begin{thm}
	Let $\De$ denote the subset of $n$-dimensional real projective space represented by vectors with non-negative entries, and let $\{N_j\}_{j \in \calj}$ denote a collection of injective $(n+1)\times(n+1)$ matrices acting projectively on $\De\to \De$. Suppose that $\calg \subset \De$ is the unique non-empty closed subset which is invariant under this action:
		$$
				\calg
			=
				\bigcup_{j \in \calj}
					N_j \cdot \calg
			\subset
				\De
			,
		$$
	and suppose that $\vol_{n-1}(\calg) = 0$, the interiors of the $N_j\cdot \De$ are mutually disjoint, and the complement $\De \setminus\bigcup_{j\in \calj} N_j \cdot \De$ is comprised by convex sets not meeting the boundary of $\De$.

	Then, indexing the connected components of $\De\setminus \calg$ by $\{\nabla\}_{i \in \cali}$, the upper box-counting dimension of $\calg$ is given by
		$$
				\overline\dim_B(\calg)
			=
				\inf
				\left\{
					s \in [n-1,n]
				\colon
					\sum_{i \in \cali}
						\vol(\nabla_i)
						\In(\nabla_i)^{s-n}
				<
					\infty
				\right\}.
		$$
\end{thm}

One can further use this to prove the following, affirming a conjecture of De Leo under these reasonably broad assumptions.
\begin{cor}
	Under the above assumptions; if in particular $\det(N_j) = \pm 1$ for each $j \in \calj$ (this can be assumed without loss of generality),
		$$
				\overline\dim_B(\calg)
			\geq
				\inf
					\left\{
						s \in [n-1,n]
					\colon
						\sum_{m=1}^\infty
						\sum_{j\in \calj^m}
							\|N_{j_1}N_{j_2}\cdots N_{j_m}\|^{-(n+1)s/n}
					<
						\infty
					\right\}.
		$$
\end{cor}


\begin{thebibliography}{2}

\bibitem{apostol}
	T. M. Apostol and M. A. Mnatsakanian.
	Figures circumscribing circles.
	\textit{Am. Math. Monthly} \textbf{111} (2014)
	853--863.

\bibitem{aspegren}
	J. Aspegren.
	On surface measures of convex bodies and generalizations of known tangential identities.
	\textit{Algebras, Groups and Geometries}
	\textbf{35} (2008) no.3, 231--242.

\bibitem{sewell}
	B. Sewell.
	\textit{A formula for the upper box-counting dimension of self-projective sets.} \verb|arXiv:2306.03047|



\end{thebibliography}
\end{document}